\documentclass[10.5pt]{article}
\usepackage{latexsym,amsfonts,euscript,amsmath}
\usepackage{xcolor}
\usepackage{amssymb}
\usepackage{stmaryrd}
\usepackage{mathrsfs}
\usepackage{leftidx}

\def\pf{\noindent{\bf Proof\quad }}
\def\pfend{\hfill{$\Box$}\bigskip}

\newtheorem{lem}{\bf Lemma}[section]

\newtheorem{pro}[lem]{\bf Proposition}
\newtheorem{thm}[lem]{\bf Theorem}

\newtheorem{cor}[lem]{\bf Corollary}

\def\la{\langle}\def\ra{\rangle}

\makeatletter 
\@addtoreset{equation}{section}
\makeatother  

\title{\bf Character codegrees, kernels, and Fitting heights of solvable groups
\thanks{Project supported by the  NSF of China (Nos. 12171058, 12301018), the NSF of Jiangsu Province (No. BK20231356) and the Natural Science Foundation of the Jiangsu Province Higher Education Institutions of China (No.
23KJB110002).}
}

\author{Guohua Qian\,\, Yu Zeng\\
{\footnotesize\small  Dept. of  Mathematics, Changshu Institute of
Technology, Changshu, Jiangsu 215500, China}\\
{\footnotesize\small E-mail: ghqian2000@163.com}, {\footnotesize\small  yuzeng2004@163.com}}

\begin{document}
\maketitle
\date{}

\vskip 1cm

\textbf{Abstract}\,\, For  an irreducible character $\chi$ of a  finite group $G$,
let $\mathrm{cod}(\chi):=|G: \ker(\chi)|/\chi(1)$ denote the codegree of $\chi$,
and let  $\mathrm{cod}(G)$ be the set of irreducible character codegrees of $G$.
In this note, we prove that if $\ker(\chi)$ is not nilpotent,
 then there exists an irreducible character $\xi$ of $G$ such that
 $\ker(\xi)<\ker(\chi)$ and $\mathrm{cod}(\xi)> \mathrm{cod}(\chi)$.
 This provides a character codegree analogue of a classical theorem of Broline and Garrison.
 As a consequence, we obtain that for a nonidentity solvable group $G$,
 its Fitting height $\ell_{\mathbf{F}}(G)$  does not exceed $|\mathrm{cod}(G)|-1$.
 Additionally,  we provide two other upper bounds for the Fitting height of a solvable group
$G$ as follows: $\ell_{\mathbf{F}}(G)\leq \frac{1}{2}(|\mathrm{cod}(G)|+2)$, and $\ell_{\mathbf{F}}(G)\leq 8\log_2(|\mathrm{cod}(G)|)+80$.


\textbf{2020 MR Subject Classification}\,\,20C15

\vskip 5cm

\pagebreak

\section{Introduction}

In this note,
$G$ always denotes a finite group,
and all characters are complex characters.
We follow Isaacs \cite{I} for standard notation and results from character theory.
For an irreducible character $\chi$ of $G$,
its codegree  is defined in \cite{Q02} as $${\rm cod}(\chi)=\frac{|G: \ker(\chi)|}{\chi(1)}.$$
Due to the natural connection between character codegrees and character degrees,
the latter can provide useful insights into the former.

A classical theorem of Broline and Garrison (\cite[Theorem 12.19]{I}) asserts that, for an irreducible character $\chi$ of $G$,
if $\ker(\chi)$ is not nilpotent, or if  $\ker(\chi)=\mathbf{F}(G)< \mathrm{Sol}(G)$,
then $G$ admits an irreducible character $\xi$  such that $\ker(\xi)<\ker(\chi)$ and $\xi(1)>\chi(1)$.
Here, $\mathbf{F}(G)$ and $\mathrm{Sol}(G)$ denote the Fitting subgroup and the solvable radical of $G$, respectively.
This theorem by Broline and Garrison has a fully dual version in terms of character codegrees,
as stated in the following result.

\begin{thm}\label{t101}  Let $\chi\in \mathrm{Irr}(G)$ and set $K=\ker(\chi)$, assume  that

{\rm (1)} either $K\nleq \mathbf{F}(G)$;

{\rm (2)} or $K=\mathbf{F}(G)<\mathrm{Sol}(G)$.

\noindent Then $G$ admits an irreducible character $\xi$ such that
$\ker(\xi)<K$ and $\mathrm{cod}(\xi)>\mathrm{cod}(\chi)$. \end{thm}

Let $\mathrm{cd}(G)$ and $\mathrm{cod}(G)$ represent  the set of irreducible character degrees and the set of irreducible character codegrees,
respectively, of $G$.
As an immediate consequence of Broline-Garrison theorem,
we know  that for a solvable group $G$,
 the Fitting height of $G$, denoted by $\ell_{\mathbf{F}}(G)$, does not exceed the cardinality of
 $\mathrm{cd}(G)$ (see, for instance, \cite[Corollary 12.21]{I}).
Analogously, Theorem \ref{t101} leads to the following corollary.

\begin{cor}\label{c102} Let $G$ be a solvable group. Then $\ell_{\mathbf{F}}(G)\leq |\mathrm{cod}(G)|-1$.\end{cor}

Note that if $G=1$, we stipulate $\ell_{\mathbf{F}}(G)=0$.

For a solvable group $G$, when
$|\mathrm{cod}(G)|$ is slightly larger, the upper bound for $\ell_{\mathbf{F}}(G)$
 provided in Corollary \ref{c102} becomes rather coarse.
 By leveraging the classification theorem for $\frac{1}{2}$-transitive permutation groups (\cite[Corollary 2]{LPS})
 and Keller's theorem for solvable linear groups (\cite[Theorem 4.5]{K}), we can obtain more precise results as follows.

\begin{pro}\label{p103} Let $G$ be a solvable group. Then $\ell_{\mathbf{F}}(G)\leq \frac{1}{2}(|\mathrm{cod}(G)|+2)$.\end{pro}

\begin{pro}\label{p104} Let $G$ be a solvable group. Then $\ell_{\mathbf{F}}(G)\leq 8\log_2(|\mathrm{cod}(G)|)+80$.\end{pro}

Examining the symmetric group $\mathsf{S}_4$, we see that the inequalities
in Corollary \ref{c102} and Proposition \ref{p103}  can both be attained as equalities.

\section{Theorem \ref{t101} and Corollary \ref{c102}}

Let $G$ be a finite group and $N$ its normal subgroup.
We identify $\chi \in \mathrm{Irr}(G/N)$ with its inflation and view $\mathrm{Irr}(G/N)$ as a subset of $\mathrm{Irr}(G)$.
We also use $\mathrm{Irr}(G|N)$ to denote the complement of $\mathrm{Irr}(G/N)$ in $\mathrm{Irr}(G)$.
Set $\mathrm{cod}(G/N)=\{ \mathrm{cod}(\chi): \chi \in \mathrm{Irr}(G/N) \}$ and $\mathrm{cod}(G|N)=\{ \mathrm{cod}(\chi):\chi \in \mathrm{Irr}(G|N) \}$.
Also, for a character $\mu$ of $G$, we use $\mathrm{Irr}(\mu)$ to denote the set of irreducible constituents of $\mu$.

We begin by presenting a useful lemma from \cite{QWW}.

\begin{lem}\label{l201} {\rm (\cite{QWW})} For an irreducible character $\chi$ of $G$, the following statements hold:

{\rm (1)} If $N$ is a $G$-invariant subgroup of $\ker (\chi)$,
then $\chi$ may be viewed as an irreducible character of $G/N$,
and the codegrees of $\chi$ in $G$ and in $G/N$ are the same;

{\rm (2)} If $M$ is subnormal in $G$ and $\psi\in \mathrm{Irr}(\chi_M)$,
then $\mathrm{cod}(\psi)$ divides $\mathrm{cod}(\chi)$.\end{lem}

The next result is fundamental and is applied in the proof of \cite[Theorem 12.19]{I}.

\begin{lem} \label{l202} Let $\chi\in \mathrm{Irr}(G)$ and $H\leq G$.
If $G=H\ker(\chi)$, then $\chi_H\in \mathrm{Irr}(H)$. \end{lem}

\pf Let $\overline{G}=G/\ker(\chi)$ and $\{ h_1,\dots ,h_t \}$ a transversal of $H\cap \ker(\chi)$ in $H$.
    Note that $\chi\in \mathrm{Irr}(\overline{G})$ and that $|G:\ker(\chi)|=|H:H\cap \ker(\chi)|=t$. Hence
	\[
		[\chi_H,\chi_H]=\frac{1}{|H|}\sum_{h\in H} \chi(h)\chi(h)^{-1}=\frac{1}{t} \sum_{i=1}^t \chi(\overline{h_i})\chi(\overline{h_i})^{-1}=\frac{1}{|\overline{G}|}	\sum_{\overline{g}\in \overline{G}} \chi(\overline{g})\chi(\overline{g})^{-1}=1,
	\]
meaning that $\chi_H\in \mathrm{Irr}(H)$.\pfend

\begin{lem} \label{l203} Let $H$ be a maximal subgroup of $G$,
let $\theta\in  \mathrm{Irr}(H)$ and $\chi\in \mathrm{Irr}(\theta^G)$.
Then $\mathrm{cod}(\chi)\geq \mathrm{cod}(\theta)$, and with equality if and only if one of the following conditions holds:

{\rm (1)} $\chi=\theta^G$ and $\ker(\chi)=\ker(\theta)$;

{\rm (2)} $\chi_H=\theta$ and $G=H \ker(\chi)$.
\end{lem}

\pf Let $K=\ker(\chi)$.
Since $H$ is a maximal subgroup of $G$,
either $K\leq H$ or $G=HK$.
Suppose that $K\leq H$.
Then $K\leq \ker(\theta)$ and $\chi(1)\leq \theta^G(1)=|G: H|\theta(1)$.
Consequently,
$$\mathrm{cod}(\chi)=\frac{|G:K|}{\chi(1)}\geq \frac{|G: \ker(\theta)|}{|G:H|\theta(1)}=\mathrm{cod}(\theta).$$
This also implies that $\mathrm{cod}(\chi)=\mathrm{cod}(\theta)$ if and only if $K=\ker(\theta)$ and $\theta^G=\chi$.
Suppose that $G=HK$. By Lemma \ref{l202}, we have $\chi_H=\theta$ and thus $\ker(\theta)=H\cap K$.
In this case, $\mathrm{cod}(\chi)=|G:K|/\chi(1)=\mathrm{cod}(\theta)$.
\pfend

\begin{pro}\label{p204} Let $H\leq G$, $\theta\in \mathrm{Irr}(H)$, and $\chi\in \mathrm{Irr}(\theta^G)$.
Then $\mathrm{cod}(\chi)\geq \mathrm{cod}(\theta)$.\end{pro}

\pf It follows by  Lemma \ref{l203} and induction.\pfend

For an irreducible character $\theta$ of a finite group $H$,
we define $\mathrm{V}(\theta):=\la h\in H|\, \theta(h)\not=0\ra$.
Clearly, $\mathrm{V}(\theta)$ is a normal subgroup of $H$.

\begin{lem}\label{l205} {\rm (\cite[Lemma 12.17]{I})}
Let $H\leq G$ and $\theta\in \mathrm{Irr}(H)$.
Suppose for each $\chi\in \mathrm{Irr}(\theta^G)$ that $\chi_H=\theta$.
Then $\mathrm{V}(\theta)\unlhd G$.\end{lem}

\begin{lem}\label{l206} {\rm (\cite[Lemma 12.18]{I})} Let $\chi\in \mathrm{Irr}(G)$ and let $N\unlhd G$ with $\ker(\chi)< N$.
Then $\ker(\chi)< N\cap \mathrm{V}(\chi)$.\end{lem}

Now, we are ready to prove Theorem \ref{t101} and then Corollary \ref{c102}.

\medskip

\noindent\textbf{Proof of Theorem \ref{t101}.}
Let $L\unlhd G$ be such that
$L=K$ if the condition (1) holds, or $L/K$ is an abelian chief factor of $G$ if the condition (2) holds.
In both cases, $L$ is not nilpotent. Hence, we can take a Sylow $p$-subgroup $P$ of $L$ such that $P$ is not normal in $L$.
Note that under condition (2), we have $L=PK$.
Applying the Frattini argument, we have in both cases that
$G=K \mathbf{N}_G(P)$.
Consequently, $G=K H$, where $H$ is a maximal subgroup of $G$ containing $\mathbf{N}_G(P)$.
Note that $\theta:=\chi_H$ is irreducible by Lemma \ref{l202}, and that $\mathrm{cod}(\chi)=\mathrm{cod}(\theta)$ by Lemma \ref{l203}.

Suppose that $\phi_H=\theta$ for all $\phi\in \mathrm{Irr}(\theta^G)$.
By Lemma \ref{l205},
$\mathrm{V}(\theta)\unlhd G$.
Observe that $$K\cap H=\ker(\theta)\leq \mathrm{V}(\theta)\leq H.$$
If $K = L$, then $P \leq K\cap H\leq L\cap \mathrm{V}(\theta)$.
Now assume that $L/K$ is an abelian chief factor of $G$.
Note that $(L\cap H)/(K\cap H)$ is also a chief factor of $H$.
Applying Lemma \ref{l206} to $\{H, \theta\}$, we obtain that $L\cap H\leq \mathrm{V}(\theta)$.
Consequently, $P\leq L\cap H\leq L\cap \mathrm{V}(\theta)$.
At this point, $P$ is always a Sylow subgroup of the $G$-invariant subgroup $L\cap \mathrm{V}(\theta)$,
and the Frattini argument shows that $G=\mathbf{N}_G(P) (L\cap \mathrm{V}(\theta))\leq H$,
which leads to a contradiction.

Therefore, there exists $\xi\in \mathrm{Irr}(\theta^G)$ such that $\xi_H\not=\theta$.
Since $$\theta^G=\chi+\xi+\cdots\not=\xi$$ and $\xi_H\not=\theta$,  we conclude by Lemma \ref{l203} that
$\mathrm{cod}(\xi)> \mathrm{cod}(\theta)$.
Therefore, $\mathrm{cod}(\xi)>\mathrm{cod}(\chi)$.
Moreover, since $\xi_H$ is not irreducible, by Lemma \ref{l202} we have $\ker(\xi)\leq H$.
Consequently, $$\ker(\xi)\leq \ker(\theta)=\ker(\chi)\cap H=K\cap H< K.$$
So, $\xi$ meets the requirements. \pfend

\noindent\textbf{Proof of Corollary \ref{c102}.} We proceed by induction on $G$.
Note that for a nonidentity group $G$, we always have
$|\mathrm{cod}(G)|\geq 2$.
Thus, we may assume  $G> \mathbf{F}(G)$.
Let $m=\max\{\mathrm{cod}(G)\}$, and let $\chi\in \mathrm{Irr}(G)$ be of codegree $m$.
By Theorem \ref{t101}, $\ker(\chi)< \mathbf{F}(G)$.
Hence, $m\not\in \mathrm{cod}(G/\mathbf{F}(G))$,
and thus $|\mathrm{cod}(G/\mathbf{F}(G))|\leq |\mathrm{cod}(G)|-1$.
Applying the inductive hypothesis to $G/\mathbf{F}(G)$,
we conclude that $\ell_{\mathbf{F}}(G/\mathbf{F}(G))\leq |\mathrm{cod}(G/\mathbf{F}(G))|-1\leq |\mathrm{cod}(G)|-2$.
Hence, $\ell_{\mathbf{F}}(G)\leq |\mathrm{cod}(G)|-1$.\pfend

\section{Propositions \ref{p103} and \ref{p104}}

Let $G$ be a finite group.
If a group $H$ is isomorphic to a subgroup of $G$, then we use the notation
$H\lesssim G$.
For a positive integer $n$ and a prime $p$, we write $n_p$ to
denote the maximal $p$-power divisor of $n$.

\begin{lem}\label{l301}
Let $p$ be a prime, and let  $G$ be a $p$-solvable group
with a unique minimal normal subgroup, say $V$. If $V=\mathbf{O}_p(G)$,
  then $\mathrm{cod}(G/V)\cap \mathrm{cod}(G|V)=\varnothing$.
  In particular, $\mathrm{cod}(G)$ is a disjoint union of $\mathrm{cod}(G/V)$ and $\mathrm{cod}(G|V)$.
\end{lem}
\pf Since $V=\mathbf{O}_{p}(G)$ is the unique minimal normal subgroup of the $p$-solvable group $G$,
 we have $G/V\lesssim\mathrm{GL}(V)$.
 It follows by \cite[Corollary 2.1]{W} that $|V|>|G/V|_p$.
 Let $\chi\in \mathrm{Irr}(G|V)$. By the unique minimal normality of $V$,
 $\chi$ is necessarily faithful.
  Hence, $\mathrm{cod}(\chi)=|G|/\chi(1)=|V|\cdot (|G/V|/\chi(1))$.
  As $V$ is abelian, $|G/V|/\chi(1)$ is a positive integer.
  Therefore, $\mathrm{cod}(\chi)_p\geq |V|>|G/V|_p$.
  Noting that $\mathrm{cod}(\omega)_p\leq |G/V|_p$ for all $\omega\in \mathrm{Irr}(G/V)$,
  we conclude that $\mathrm{cod}(G/V)\cap \mathrm{cod}(G|V)=\varnothing$.
  \pfend

Let $V$ be a $d$-dimensional vector space over the prime field $\mathbb{F}_p$.
As in \cite{manzwolfbook}, we denote by $\Gamma(p^{d})$ the semilinear group of $V$, i.e. (identifying $V$ with $\mathbb{F}_{p^d}$)
\[
	\Gamma(p^{d})=\{ x\mapsto ax^\sigma :  x\in \mathbb{F}_{p^d}, a\in \mathbb{F}_{p^d}^\times, \sigma\in \mathrm{Gal}(\mathbb{F}_{p^d}/\mathbb{F}_p)\}.
\]
It is noteworthy that $\Gamma(p^{d})$ is a metacyclic group.

To prove Proposition \ref{p103}, we need the next result: the classification of solvable $\frac{1}{2}$-transitive permutation groups.

\begin{thm}\label{t302}
  Let $V$ be an elementary abelian group of order $p^d$ for a prime $p$,
  and $G$ a solvable subgroup of $\mathrm{GL}(V)=\mathrm{GL}(d,p)$.
  If $G$ acts $\frac{1}{2}$-transitively on $V^\sharp:=V-\{ 0 \}$,
  meaning that $|\mathbf{C}_{G}(v)|$ is a constant for each $v\in V^\sharp$,
  then one of the following holds:

{\rm (1)} $G$ is transitive on $V^\sharp$;

{\rm (2)} $G\leq \Gamma(p^d)$;

{\rm (3)} $G$ acts fixed-point-freely on $V$;

{\rm (4)} $G$ is isomorphic to the subgroup of $\mathrm{GL}(2,p^{d/2})$ of order $4(p^{d/2}-1)$ consisting of all monomial matrices of determinant $\pm 1$ where $p>2$;

{\rm (5)} $G$ is solvable and $p^d=3^2,5^2,7^2,11^2,17^2$ or $3^4$.
  \end{thm}
\pf  This is a partial result of \cite[Corollary 2]{LPS}. \pfend

The next theorem, which is known as the classification of solvable doubly transitive permutation groups, describes $VG$ if the case (1) of Theorem \ref{t302} holds.
In the following, we use the symbol $\mathsf{PrimitiveGroup}(m,i)$ for the $i$-th primitive group of degree $m$ in the $\mathsf{GAP}$ library of primitive groups (\cite{gap}).

\begin{thm}\label{t303}
  Let $V$ be an elementary abelian group of order $p^d$ for a prime $p$,
  and $G$ a solvable subgroup of $\mathrm{GL}(V)=\mathrm{GL}(d,p)$.
  If $G$ acts transitively on $V^\sharp$,
  then either $G\leq\Gamma(p^d)$, or $VG=\mathsf{PrimitiveGroup}(p^d,i)$ with
  \[
    \begin{aligned}
      (p^d, i) \in & \{(3^2, 6),(3^2, 7),(5^2, 15),(5^2, 18),(5^2, 19),(7^2, 25),(7^2, 29),\\
        &  (11^2, 39),(11^2, 42),(23^2, 59),(3^4, 71),(3^4, 90),(3^4, 99)\}
    \end{aligned}
    \]
\end{thm}
\pf
This is \cite[Hauptsatz]{H} combined with group identification via $\mathsf{GAP}$ \cite{gap}.
\pfend

In the following,
we use the symbol $\mathsf{Q}_{2^n}$ for the generalized quaternion group of order $2^n$ where $n\geq 3$, $\mathsf{S}_m$ for the symmetric group of degree $m$, $\mathsf{C}_{m}$ for the cyclic group of order $m$.
If $G$ is a $p$-solvable group, we use $\ell_p(G)$ to denote the $p$-length of $G$.

\begin{lem}\label{l303}
  Let $V$ be an elementary abelian group of order $p^d$ for a prime $p$,
  and $G$ a solvable subgroup of $\mathrm{GL}(V)=\mathrm{GL}(d,p)$.
  If $G$ acts $\frac{1}{2}$-transitively on $V^\sharp$,
  then one of the following holds:

{\rm (1)} $\ell_{\mathbf{F}}(G)\leq 2$;

{\rm (2)} $\ell_{\mathbf{F}}(G)=3$ and $|\mathrm{cod}(G)|\geq 6$;

{\rm (3)} $\ell_{\mathbf{F}}(G)=4$ and $|\mathrm{cod}(G)|\geq 8$.
\end{lem}
\pf
We proceed the proof by checking the groups in Theorem \ref{t302} case by case.

Assume that either $G\leq \Gamma(p^d)$, or $G$ is isomorphic to the subgroup of $\mathrm{GL}(2,p^{d/2})$ of order $4(p^{d/2}-1)$ consisting of all monomial matrices of determinant $\pm 1$ where $p>2$.
Then $G$ is metabelian. In particular, $\ell_{\mathbf{F}}(G)\leq 2$.

Assume that the solvable group $G$ acts fixed-point-freely on $V$.
If $\ell_{\mathbf{F}}(G)\leq 2$, then we are done.
We now consider the case when $\ell_{\mathbf{F}}(G)\geq 3$.  Let $P\in \mathrm{Syl}_2(G)$ and $R\in \mathrm{Syl}_2(\mathbf{F}(G))$.
Since $G/\mathbf{F}(G)$ is not nilpotent,  $\mathbf{F}(G)$ cannot be cyclic.
Therefore, $\mathbf{F}(G)= R\times C$, where $R\cong \mathsf{Q}_{2^n}$ and $C \in \mathrm{Hall}_{2'}(\mathbf{F}(G))$ is cyclic.
Note that
$$G/\mathbf{F}(G)\lesssim \mathrm{Out}(\mathbf{F}(G))\cong \mathrm{Out}(R)\times \mathrm{Out}(C),$$
and that $\mathrm{Out}(R)\cong \mathrm{Out}(\mathsf{Q}_{2^{n}})$ is a $2$-group except when $n=3$ (see, for instance, \cite[Section 5.3, Lemma 5.3.3 and Excercise 7]{kurzweil}).
Consequently, $R\cong \mathsf{Q}_8$ and $\mathsf{S}_3\lesssim G/\mathbf{F}(G)\lesssim \mathsf{S}_3\times A$, where $A$ is an abelian group. 
In particular, this implies $\ell_{\mathbf{F}}(G)=3$ and $\ell_2(G)=2$.
Observe that both $P$ and $R$ are generalized quaternion groups. 
It follows that $P/R$ is cyclic,
resulting in $|P: R|=|G: \mathbf{F}(G)|_2=2$.  
Thus, $P\cong \mathsf{Q}_{16}$.
Let $D\unlhd G$ be maximal such that $R\cap D=1$.
Note that $D$ has odd order, because otherwise it would imply $P\cong R\times \mathsf{C}_2$.
This leads to $\ell_2(G/D)=\ell_2(G)=2$, so $\ell_{\mathbf{F}}(G/D)=3$.
Now it is straightforward to see that
\begin{center}
  $\mathbf{F}(G/D)=RD/D
\cong R\cong \mathsf{Q}_{8}$ and $G/RD\cong \mathsf{S}_3$.
\end{center}
Now, consider $\overline{G}=G/D$.
We have that $\ell_{\mathbf{F}}(\overline{G})=\ell_{\mathbf{F}}(G/D)=3$ and $|\overline{G}|=48$. 
Note that there are exactly two non-isomorphic groups of order 48 with Fitting subgroups isomorphic to $\mathsf{Q}_8$, 
namely the general linear group $\mathrm{GL}(2,3)$ and the conformal special unitary group $\mathrm{CSU}(2,3)$.
Since Sylow $2$-subgroups of $\mathrm{GL}(2,3)$ are not isomorphic to $\mathsf{Q}_{16}$, we deduce that $\overline{G}\cong \mathrm{CSU}(2,3)$. 
Calculating the character codegrees of $\overline{G}$, we have $\mathrm{cod}(\overline{G})=\{1,2,3,8,12,24\}$.
Therefore, $|\mathrm{cod}(G)|\geq |\mathrm{cod}(\overline{G})|=6$.

Assume that the solvable group $G$ acts transitively on $V^\sharp$ and that $G$ is not isomorphic to a subgroup of $\Gamma(p^{d})$.
Theorem \ref{t303}
yields that
$VG=\mathsf{PrimitiveGroup}(p^d,i)$, where
\[
\begin{aligned}
  (p^d, i) \in & \{(3^2, 6),(3^2, 7),(5^2, 15),(5^2, 18),(5^2, 19),(7^2, 25),(7^2, 29),\\
    &  (11^2, 39),(11^2, 42),(23^2, 59),(3^4, 71),(3^4, 90),(3^4, 99)\}
\end{aligned}
\]
Checking these groups via $\mathsf{GAP}$ \cite{gap}, we conclude that either $\ell_{\mathbf{F}}(G)\leq 2$, or $\ell_{\mathbf{F}}(G)=3$ and $|\mathrm{cod}(G)|\geq 6$.

Assume finally that $G$ is solvable and $p^d\in \{3^2,5^2,7^2,11^2,17^2,3^4\}$.
Checking these groups via $\mathsf{GAP}$ \cite{gap}, we conclude that one of (1), (2) or (3) holds.
\pfend

Now, we are ready to prove Proposition \ref{p103}.

\medskip

\noindent\textbf{Proof of Proposition \ref{p103}.}   We proceed by induction on $|G|$.
  By induction, we may assume that $\Phi(G)=1$ and that $V:=\mathbf{F}(G)$ is a unique minimal normal subgroup of $G$.
  Therefore, there exists a maximal subgroup $H$ of $G$ such that  $G=H \ltimes V$.
  By induction, $\ell_{\mathbf{F}}(G/V)\leq \frac{1}{2}(|\mathrm{cod}(G/V)|+2)$.
  Thus,
  \[
    \ell_{\mathbf{F}}(G)=1+\ell_{\mathbf{F}}(G/V)\leq \frac{1}{2}(|\mathrm{cod}(G/V)|+4).
  \]
  Note that by Lemma \ref{l301}, $\mathrm{cod}(G/V)\cap \mathrm{cod}(G|V)=\varnothing$.
  If $|\mathrm{cod}(G|V)|\geq 2$, then
  $$\ell_{\mathbf{F}}(G)\leq \frac{1}{2}(|\mathrm{cod}(G/V)|+4)\leq \frac{1}{2}(|\mathrm{cod}(G/V)|+|\mathrm{cod}(G|V)|+2)=\frac{1}{2}(|\mathrm{cod}(G)|+2),$$
  and we are done.
  
  We now assume that $|\mathrm{cod}(G|V)|=1$.
   Since $V$ is the unique minimal normal subgroup of $G$,
   all $\chi \in \mathrm{Irr}(G|V)$ are faithful and therefore have the same degree.
   For a character $\lambda \in \mathrm{Irr}(V)$, as $G=H \ltimes V$,
   $\lambda$ extends to $\mathrm{I}_{G}(\lambda)$.
   By Clifford's theorem, $|H:\mathrm{I}_{H}(\lambda)|$ is a constant for each nontrivial $\lambda \in \mathrm{Irr}(V)$.
  Consequently, $H$ acts $\frac{1}{2}$-transitively on $\mathrm{Irr}(V)-\{ 1_V \}$.
 
  Now, assume that $\ell_{\mathbf{F}}(G)\leq 3$. Then $\ell_{\mathbf{F}}(G)\leq \frac{1}{2}(\ell_{\mathbf{F}}(G)+3)$,
  and  Corollary \ref{c102} implies that $ \ell_{\mathbf{F}}(G)\leq \frac{1}{2}(|\mathrm{cod}(G)|-1+3)=\frac{1}{2}(|\mathrm{cod}(G)|+2)$,
  as desired.

  Assume next that $\ell_{\mathbf{F}}(G)\geq 4$.
  By Lemma \ref{l303},
  we have that either $\ell_{\mathbf{F}}(G/V)=3$ and $|\mathrm{cod}(G/V)|\geq 6$,
  or $\ell_{\mathbf{F}}(G/V)=4$ and $|\mathrm{cod}(G/V)|\geq 8$.
 Since by Lemma \ref{l301}
 $$|\mathrm{cod}(G)|=|\mathrm{cod}(G/V)|+|\mathrm{cod}(G|V)|=|\mathrm{cod}(G/V)|+1,$$ 
  we easily conclude that $ \ell_{\mathbf{F}}(G)\leq \frac{1}{2}(|\mathrm{cod}(G)|+2)$.
\pfend

Let $G$ acts on a vector space $V$,
and denote by 
 $m(G, V)$ the number of different nontrivial orbit sizes induced by the action of $G$ on $V$.
Specifically,
$m(G, V)=|\{|v^G|:\, 0\not=v\in V\}|$, where $v^G=\{v^g: g\in G\}$.
It follows from \cite[Theorem 4.5]{K} that
if a solvable group $G$ acts faithfully and irreducibly on a finite vector space $V$, then
\begin{equation}\label{eq301} \mathrm{dl}(G/\mathbf{F}(G))\leq 8\,\log_2 (m(G, V))+78,\end{equation}
where $\mathrm{dl}(G/\mathbf{F}(G))$ denotes the derived length of $G/\mathbf{F}(G)$.
By applying Keller's result (\ref{eq301}),
we can  derive a logarithmic bound for the Fitting height  $\ell_{\mathbf{F}}(G)$ of solvable groups $G$,
expressed in terms of $|\mathrm{cod}(G)|$.

\begin{pro} Let $G$ be a solvable group.
Then $\mathrm{dl}(G/\mathbf{F}_2(G))\leq 8\,\log_2 (|\mathrm{cod}(G)|)+78$,
and in particular, $\ell_{\mathbf{F}}(G)\leq 8\log_2(|\mathrm{cod}(G)|)+80$.\end{pro}
\pf Note that $\mathrm{cod}(G/N)\subseteq \mathrm{cod}(G)$ whenever $N\unlhd G$.
Assume that $G$ admits different minimal normal subgroups $N_a$ and $N_b$.
Let $U_a/N_a=\mathbf{F}_2(G/N_a)$, $U_b/N_b=\mathbf{F}_2(G/N_b)$, 
where $\mathbf{F}_2(H)$ denotes the preimage of $\mathbf{F}(H/\mathbf{F}(H))$ in $H$.
Then $\mathbf{F}_2(G)=U_a\cap U_b$, and $\mathrm{dl}(G/U_c)\leq 8\,\log_2 (|\mathrm{cod}(G)|)+78$ by induction, where $c\in \{a, b\}$.
This implies that 
$$\mathrm{dl}(G/\mathbf{F}_2(G))=\max\{\mathrm{dl}(G/U_a),\mathrm{dl}(G/U_b)\} \leq 8\,\log_2 (|\mathrm{cod}(G)|)+78,$$
and we are done. Hence, we may assume that $G$ admits a unique minimal normal subgroup, say $N$.
Similarly, we may assume by induction that $\Phi(G)=1$.
Consequently, there exists a maximal subgroup $H$ of $G$ such that $G=H\ltimes N$.

Obviously, $H$ acts faithfully and irreducibly on $N$ and $V:=\mathrm{Irr}(N)$.
For each $s\in m(H,V)$,
let $1_N\not=\lambda_s\in V$ be such that  $|H: \mathrm{I}_H(\lambda_s)|=s$.
Since $\lambda_s$ is extendible to $\mathrm{I}_G(\lambda_s)$, we can choose $\chi_s\in \mathrm{Irr}(\lambda_s^G)$ such that
$\chi_s(1)=s$. Note that $\chi_s$ is faithful, and it follows that
$|G|/s\in \mathrm{cod}(G)$. In particular,
$|m(H, V)|\leq |\mathrm{cod}(G)|$.
Applying inequality (\ref{eq301}) and observing that $\mathbf{F}_2(G)=\mathbf{F}(H)N$,
we conclude that
$\mathrm{dl}(G/\mathbf{F}_2(G))=\mathrm{dl}(H/\mathbf{F}(H))\leq  8\,\log_2 (|\mathrm{cod}(G)|)+78$, as desired.\pfend

\bigskip

\noindent\textbf{Data availability statements} Data sharing not applicable to this article as no
 datasets were generated or analysed during the current study.

 \bigskip

\noindent\textbf{Declaration}\\ 
The authors declare that they have no conflict of interest.
The authors whose names appear on the submission have contributed sufficiently to the scientific work and therefore share collective responsibility and accountability for the results.\\



\end{document}